\documentclass[12pt]{amsart}
\usepackage{latexsym,amsthm,amssymb,amscd}

\newcommand{\NN}{\Bbb N}

\newcommand{\ZZ}{\Bbb Z}
\newcommand{\mm}{\mathfrak m}
\newcommand{\nn}{\mathfrak n}

\newcommand{\reg}{{\rm reg}}

\newcommand{\indeg}{{\rm indeg}}
\newcommand{\ini}{{\rm in}}

\newcommand{\rate}{{\rm rate}}

\newcommand{\ch}{{\rm char}}

\newcommand{\Tor}{{\rm Tor}}

\newcommand{\id}{{\rm id}}

\newcommand{\Gin}{{\rm Gin}}
\newcommand{\union}{\cup}

\def\pnt{{\raise 0.5mm\hbox{\large\bf.}}}

\newtheorem{Theorem}{\bf Theorem}[section]
\newtheorem{Lemma}[Theorem]{\bf Lemma}
\newtheorem{Corollary}[Theorem]{\bf Corollary}
\newtheorem{Proposition}[Theorem]{\bf Proposition}
\newtheorem{Remark}[Theorem]{\bf Remark}

\newtheorem{Example and Notation}[Theorem]{\bf Example and Notation}

\newtheorem{Theorem and Definition}[Theorem]{\bf Theorem and Definition}

\let\epsilon=\varepsilon
\let\phi=\varphi
\let\kappa=\varkappa
\let\to=\rightarrow

\title{Subalgebras of bigraded Koszul Algebras}
\author{Stefan Blum}

\begin{document}

\begin{abstract}
We show that diagonal subalgebras and generalized Veronese
subrings of a bigraded Koszul algebra are Koszul.
We give upper bounds 
for the regularity of sidediagonal and relative Veronese modules 
and apply the results to symmetric algebras and Rees rings. 
\end{abstract}

\maketitle

\section*{Introduction}
In this article we study standard bigraded algebras. Let $K$ denote a field
and let  $R=S/J$ where $S=K[x_1,\ldots,x_n,y_1,\ldots,y_m]$ is a polynomial ring 
with standard bigrading $\deg(x_i)=(1,0)$ and $\deg(y_j)=(0,1)$ and $J \subset S$ a 
bihomogeneous ideal. For such an algebra we consider two kinds of subalgebras. 
Let $a,b \geq 0$ be two integers with $(a,b) \neq (0,0)$. Then the $(a,b)$-diagonal 
subalgebra is the positively graded subring 
$$R_\Delta = \bigoplus_{i \geq 0} R_{(ia,ib)}$$ 
where $R_{(i,j)}$ denotes the $(i,j)^{th}$ bigraded component of $R$. Moreover, for two
integers $a,b \geq 0$ such that $(a,b) \neq (0,0)$, a generalized bigraded Veronese subring 
of $R$ is defined by 
$$R_{\Tilde{\Delta}} = \bigoplus_{i,j \geq 0} R_{(ia,jb)}.$$ 
Recall that a positively graded $K$-algebra $A$ is called Koszul if the
residue class field $K$, considered as a trivial $A$-module, has linear $A$-free
resolution. During the last thirty years Koszul algebras have been studied in various 
contexts. A good survey is given by Fr\"oberg in \cite{F}.
       
In the last years, diagonal subalgebras have been intensively studied because 
they naturally appear in Rees algebras and symmetric algebras. In \cite{CHTV} Conca, 
Herzog, Trung and Valla discuss many algebraic properties.
In particular, they prove that for an arbitrary bigraded algebra $R$, the diagonals
$R_\Delta$ are Koszul provided one chooses $a$ and $b$ large enough. In this article
they ask two questions of which one is positively answered by Aramova, Crona and De Negri 
in \cite{ACD} who showed that the defining ideal of $R_\Delta$ has quadratic Gr\"obner 
basis for $a,b \gg 0$. It is a well-known fact that this is a stronger property than 
Koszulness. The main result of this article is the positive answer to the second question 
posed in \cite{CHTV}: Suppose $R$ is a Koszul algebra, then all diagonal subalgebras 
$R_\Delta$ are Koszul. Moreover, we prove that all generalized Veronese subrings 
$R_{\Tilde{\Delta}}$ inherit the Koszul property. The algebras $R_{\Tilde{\Delta}}$    
appear in \cite{R} where R\"omer studies homological properties of bigraded algebras.
  
Note that the Castelnouvo-Mumford regularity (see Section 1 for the definition) measures 
the complexity of the minimal free resolution of a finitely generated $R$-module. 
Provided $R$ is a Koszul algebra all finitely generated modules have finite regularity 
over $R$ (see \cite{AE}).

For a finitely generated bigraded $R$-module and two integers 
$c,d \geq 0$ we define a sidediagonal module $M^{(c,d)}_\Delta$ as the $R_\Delta$-module with graded 
components $(M^{(c,d)}_\Delta)_{i}=M_{(ia+c,ib+d)}$. Similarly one defines a bigraded 
relative Veronese  module $M^{(c,d)}_{\Tilde{\Delta}}$. 
Provided $R$ is Koszul, we get upper bounds for the regularity of these modules. 

In particular, if the initial degree of $M$ (see Section 1 for the definition) equals $0$  
and if we choose $c,d$ such that $R^{(c,d)}_\Delta$ respectively  $R^{(c,d)}_{\Tilde{\Delta}}$ 
are generated in degree $0$, then $\reg_{R_\Delta} M^{(c,d)}_\Delta \leq 1$ and
$\reg_{R_{\Tilde{\Delta}}} M^{(c,d)}_{\Tilde{\Delta}} \leq 2$ for all
$\Delta$ respectively $\Tilde{\Delta}$ with $a,b \geq \reg_R M$.
For the proof we use similar techniques as Aramova, Barcanescu, and Herzog in \cite{ABH} 
where they give upper bounds for rates of modules over arbitrary Veronese algebras.   
Notice that our results also hold with similar proofs if one considers multigraded 
$K$-algebras and the corresponding multigraded subalgebras. 

This paper is structured in the following way. In the first section we recall definitions 
and introduce notation. 

In Section $2$ we prove the main result and get the upper bounds for the regularities 
mentioned above.
 
In the third section we discuss some applications of the main result. Let $A$ be a
positively graded $K$-algebra and $M$ a finitely generated $A$-module.
Provided the symmetric algebra $S(M)$ is Koszul, we show that all symmetric powers
of $M$ have linear resolutions. In the specific case that $M =\mm$ is the graded
maximal ideal of $A$, we have a necessary condition for $S(\mm)$ to 
be Koszul, that is when the defining ideal of $A$ has a $2$-linear resolution.
Under the weaker assumption that $A$ is Koszul, we obtain that all
symmetric powers of $\mm$ have a linear $A$-resolution.
Let $I \subset A$ be an homogeneous ideal generated in one degree. If the Rees ring $R(I)$ 
is Koszul, then all powers of $I$ have linear $A$-resolutions.

Moreover, we recover some well-known results about graded Koszul algebras which
Backelin and Fr\"oberg first proved in \cite{BF} saying that the Koszul property is 
preserved for tensor products over $K$, Segre products and Veronese subrings.
  
In the last section we interpret the main result for semigroup rings. The
Koszul property for these rings corresponds to the Cohen-Macaulay property of
certain divisor posets (see \cite{HRW} and \cite{PRS}). Therefore, we obtain 
that Cohen-Macaulayness for these divisor posets is preserved under
taking diagonals and generalized Veronese subrings.

The author is grateful to Prof.\ Herzog for several inspiring discussions on the
subject of this article.  

\section{Notation}
Let $K$ be a field and $S=K[x_1,\ldots,x_n,y_1,\ldots,y_m]$ be the polynomial ring with 
standard bigrading $\deg(x_i)=(1,0)$ and $\deg(y_j)=(0,1)$. Throughout this paper $R$ 
always denotes a bigraded $K$-algebra of the form $R=S/J$ where $J$ is a bihomogeneous 
ideal of $S$. We recall that for two integers $a,b \geq 0$ with $(a,b)\neq (0,0)$
the $(a,b)$-diagonal is the set $\Delta =\{(sa,sb) \colon s \in \ZZ \}$ of $\ZZ^2$. 
As in \cite{H} the diagonal subalgebra of $R$ along $\Delta$ is defined as the 
positively graded algebra  
$$R_\Delta=\bigoplus_{i \geq 0} R_{(ia,ib)}$$
where $R_{(i,j)}$ denotes the $(i,j)^{th}$ bigraded component of $R$. The algebra $R_\Delta$
is generated by the residue classes of all monomials which have degree $(a,b)$ in $S$.  
Therefore, $R_\Delta$ is standard graded. According to \cite{R} we define for two integers 
$a,b \geq 0$ with $(a,b)\neq (0,0)$ the bigraded generalized Veronese subring of $R$ along 
the set $\tilde{\Delta} =\{(sa,tb) \colon s,t \in \ZZ \}$ by 
$$R_{\tilde{\Delta}}=\bigoplus_{i,j \geq 0} R_{(ia,jb)}.$$
Here the bigraded components are $(R_{\tilde{\Delta}})_{(i,j)}=R_{(ia,jb)}$. 
The algebra $R_{\tilde{\Delta}}$is generated by the residue classes of all monomials
which have degree $(a,0)$ or $(0,b)$ in $S$. Thus, $R_{\tilde{\Delta}}$ is a standard
bigraded algebra. Note that $R = R_{\Tilde{\Delta}}$ 
for $(a,b)=(1,1)$. In the case that $n=0$ or $m=0$, the algebra $R$ is simply standard
graded and the subrings $R_\Delta$ resp.\ $R_{\Tilde{\Delta}}$ are the well-known Veronese
subrings of $R$. Observe also that the $(1,1)$-diagonal of $R_{\Tilde{\Delta}}$ equals 
$R_{\Delta}$.
 
Let $M$ be a finitely generated, bigraded $R$-module. For two integers $c,d \geq 0$
we define $M^{(c,d)}_{\Delta}$ to be the finitely generated, $\ZZ$-graded 
$R_{\Delta}$-module  with components $(M^{(c,d)}_{\Delta})_i=M_{(ia+c,ib+d)}$.
For $(c,d)=(0,0)$ we simply use $M_\Delta$ instead of $M^{(0,0)}_\Delta$. We call  
$M^{(c,d)}_{\Delta}$ the \it$(c,d)$-sidediagonal module \rm of $M$. Similarly, we write 
$M^{(c,d)}_{\Tilde{\Delta}}$ for the bigraded $R_{\Tilde{\Delta}}$-module with components 
$(M^{(c,d)}_{\Tilde{\Delta}})_{(i,j)}=M_{(ia+c,jb+d)}$
and call it the \it relative $(c,d)$-Veronese module \rm of $M$. If $n=0$ or $m=0$, then
these modules coincide with the relative Veronese modules defined in \cite{ABH}. 
We need two index sets 
$$\mathcal{I}(a,b) = \{ (c,d) \in \NN^2 \colon  c < a \text{ or } d < b\} \quad\text{ and}$$
$$ 
\Tilde{\mathcal{I}}(a,b) = \begin{cases}
                              \{ (c,d) \in \NN^2 \colon  c < a \text{ and } d < b\} &
							  \text{ if } a,b \geq 1, \cr
                              \{ (c,0) \in \NN^2 \colon  c < a \}  &
							  \text{ if } a \geq 1 \text{ and } b = 0, \cr
                              \{ (0,d) \in \NN^2 \colon  d < b \} &
							  \text{ if } a = 0 \text{ and } b \geq 1. \cr
							\end{cases}  
$$
Note that the index set $\mathcal{I}(a,b)$ is infinite while 
$\Tilde{\mathcal{I}}(a,b)$ is a finite set.
For $(c,d) \in \mathcal{I}(a,b)$ the module $R^{(c,d)}_{\Delta}$  is generated in 
degree $0$ and, for arbitrary $c,d \geq 0$, it is 
$R^{(c,d)}_{\Delta} = R^{(c',d')}_{\Delta}(-l)$ with some integer $l \geq 0$ and
some $(c',d') \in \mathcal{I}(a,b)$. 
An analogous fact holds for the modules $R^{(c,d)}_{\Tilde{\Delta}}$.
We have the decomposition
$$R=\bigoplus_{(c,d) \in \mathcal{I}(a,b)} R^{(c,d)}_{\Delta}.$$ 
Analogously, if $a,b \geq 1$, then $R$ is the finite direct sum of the 
$R^{(c,d)}_{\Tilde{\Delta}}$ with $(c,d) \in \Tilde{\mathcal{I}}(a,b)$.

The map $M \mapsto M^{(c,d)}_{\Delta}$ (resp.\ $M \mapsto M^{(c,d)}_{\Tilde{\Delta}}$) 
defines an exact functor from the category of bigraded finitely generated $R$-modules 
to the category of $\ZZ$-graded finitely generated $R_\Delta$-modules 
(resp.\ bigraded $R_{\Tilde{\Delta}}$-modules).
In particular, consider a bigraded free resolution 
$$ F_\pnt: \quad \ldots \to F_i \to \ldots \to F_1 \to F_0 \to M \to 0$$
where every free module $F_i$ decomposes into a finite direct 
$F_i = \bigoplus_{p,q} R(-p,-q)^{b_{i,(p,q)}}$. Here, $R(-p,-q)$ denotes the bigraded 
$R$-module with components $R(-p,-q)_{(i,j)}= R_{(i-q,j-q)}$. Then we get an exact complex
of $R_\Delta$-modules 
$$ (F_\pnt)^{(c,d)}_{\Delta}: \quad \ldots \to (F_i)^{(c,d)}_\Delta \to \ldots \to 
(F_1)^{(c,d)}_\Delta  \to (F_0)^{(c,d)}_\Delta \to M^{(c,d)}_{\Delta} \to 0$$
with $ (F_i)^{(c,d)}_\Delta = \bigoplus_{p,q} 
( R(-p,-q))^{(c,d)}_\Delta )^{b_{i,(p,q)}}$.
Analogous statements are true for the functor $-^{(c,d)}_{\Tilde{\Delta}}$. 
It will be important for the main result to write every module 
$R(-p,-q))^{(c,d)}_\Delta$ as a shifted sidediagonal module of the form 
$R^{(c',d')}_\Delta$ for some $(c',d') \in \mathcal{I}(a,b)$. 
For a reel number $\alpha$ we use $\lceil \alpha \rceil$ for the smallest integer 
$z$ such that $z \geq \alpha$. We observe the following.

\begin{Remark}\label{shifts}
\rm Let $\Delta$ be the $(a,b)$-diagonal. For $z \in \ZZ$ let 
$\alpha(z) \in \{0,\ldots,a-1\}$ be the integer such that $\alpha(z) \equiv z \mod a$
and $\beta(z) \in \{0,\ldots,b-1\}$ with $\beta(z) \equiv z \mod b$.
\begin{enumerate}
\item
\begin{enumerate}
\item
Let $a>0$, $b=0$ and $(c,d) \in \mathcal{I}(a,b)$. Then
$$R(-p,-q)^{(c,d)}_{\Delta}= \begin{cases} 0, & \text{if } q > d, \cr
                                           R^{(\alpha(c-p),d-q)}_\Delta (-l),
										    & \text{if } q \leq d, \cr
						      \end{cases} $$
where $l=\max\{0,\lceil \frac{p-c}{a}\rceil\}$.
\item
Let $a=0$, $b>0$ and $(c,d) \in \mathcal{I}(a,b)$. Then
$$R(-p,-q)^{(c,d)}_{\Delta}= \begin{cases} 0, & \text{if } p > c, \cr
                                           R^{(c-p,\beta(d-q))}_\Delta (-l),
										   & \text{if } p \leq c, \cr 
						      \end{cases} $$
where $l=\max\{0,\lceil \frac{q-d}{b}\rceil\}$.
\item
Let $a,b \geq 1$ and  $(c,d) \in \mathcal{I}(a,b)$. Then
$$R(-p,-q)^{(c,d)}_{\Delta}= R^{(la+c-p,lb+d-q)}_{\Delta}(-l)$$
where $l=\max\{0,\lceil \frac{p-c}{a}\rceil,\lceil \frac{q-d}{b}\rceil\}$.
\end{enumerate}
\item 
\begin{enumerate}
\item
Let $a>0$, $b=0$ and $(c,0) \in \Tilde{\mathcal{I}}(a,b)$. Then
$$R(-p,-q)^{(c,0)}_{\Tilde{\Delta}}= \begin{cases} 0, & \text{if } q > 0, \cr
                                           R^{(\alpha(c-p),0)}_{\Tilde{\Delta}} (-k,0),
										    & \text{if } q = 0, \cr
						      \end{cases} $$
where $k=\max\{0,\lceil \frac{p-c}{a}\rceil\}$.
\item
Let $a=0$, $b>0$ and $(0,d) \in \Tilde{\mathcal{I}}(a,b)$. Then
$$R(-p,-q)^{(0,d)}_{\Tilde{\Delta}}= \begin{cases} 0, & \text{if } p > 0, \cr
                                           R^{(0,\beta(d-q))}_{\Tilde{\Delta}} (0,-l),
										   & \text{if } p = 0,  \cr 
						      \end{cases} $$
where $l=\max\{0,\lceil \frac{q-d}{b}\rceil\}$.
\item
Let $a,b \geq 1$ and $(c,d) \in \Tilde{\mathcal{I}}(a,b)$. Then
$$R(-p,-q)^{(c,d)}_{\Tilde{\Delta}}= R^{(\alpha(c-p),\beta(d-q))}_{\Tilde{\Delta}} (-k,-l)$$
where  $k=\max\{0,\lceil \frac{p-c}{a}\rceil\}$ and $l=\max\{0,\lceil \frac{q-d}{b}\rceil\}$ 
\end{enumerate}
\end{enumerate}
\end{Remark}

We recall some well-known definitions. For a bigraded, finitely generated $R$-module $M$ 
each $\Tor^R_i(M,K)$-group is naturally bigraded and the bigraded Poincar\'e series is
given by 
$$P^R_M(s,t,z) =\sum_{i,j,k} \dim_K \Tor^R_i(M,K)_{(j,k)} s^j t^k z^i.$$
Let $A$ be a positively graded $K$-algebra and $N$ a finitely generated $A$-module. We set
$$t_s(N)=\sup\{j \colon \Tor^A_s(N,K)_j \neq 0\}$$
with $t_s(N) = -\infty$ if $\Tor^A_s(N,K)=0$. Recall that the Castelnuovo-Mumford
regularity is defined as
$$\reg_A N = \sup\{t_i(N) -i \colon i \geq 0\}.$$ 
The initial degree $\indeg(N)$ of $N$ is the minimum of the $i$ such that $N_i \neq 0$.
Note that $M$ is said to have an $i$-linear resolution if $\reg_A M = \indeg(M) =i$. 
By definition,  $A$ is Koszul if and only if $\reg_A K =0$. 
Every bigraded $K$-algebra $R$ is also naturally $\NN$-graded with $i^{th}$ component
$R_i = \bigoplus_{j+k=i} R_{(j,k)}.$ Similarly, every bigraded $R$-module $M$ can be 
considered as $\ZZ$-graded. We say that $M$ has a bigraded $a$-linear resolution if 
$\Tor^R_i(M,K)_{(j,k)} = 0$ for all $i \geq 0$ and all $j + k \neq i + a$. 
\section{Main result}
In this section we prove the main result of this article.

\begin{Theorem}\label{main}
\rm If $R$ is a Koszul algebra, then every diagonal subalgebra $R_\Delta$ and 
every generalized Veronese subring $R_{\tilde{\Delta}}$is a Koszul
algebra.
\end{Theorem}

For the proof we need several lemmata. Let  
$\nn_x = (x_1,\ldots,x_n) \subset R$ resp.\ $\nn_y =(y_1,\ldots,y_m) \subset R$ be the 
ideal generated by the residue classes of the $x_i$ resp.\ the $y_j$.

\begin{Lemma}\label{Koszul_linear}
\rm If $R$ is Koszul, then the ideals $\nn_x$ and $\nn_y$ have bigraded $1$-linear 
$R$-resolutions.
\end{Lemma}

\begin{proof}
By symmetry, it is enough to show that $\nn_y$ has a bigraded linear resolution.
The residue class field $K$ has a $0$-linear minimal free $R$-resolution $F_\pnt$
because $R$ is Koszul. Let $\Delta$ be the $(1,0)$-diagonal. 
Applying the functor $-_{\Tilde{\Delta}}$ we get the exact complex 
$(F_\pnt)_{\Tilde{\Delta}} \to K \to 0$. 
By Remark \ref{shifts}(b), the $i^{th}$ module $(F_i)_{\Tilde{\Delta}}$ 
is a direct sum of copies of $R_{\Tilde{\Delta}}$ shifted by $(-i,0)$. 
Thus $R_{\Tilde{\Delta}}$ is a standard bigraded Koszul algebra.

Let $p:R \to R_{\Tilde{\Delta}}$ be the projection map and $i: R_{\Tilde{\Delta}} \to R$ 
the inclusion. Note that both maps $p$ and $i$ are bigraded homomorphisms.
Since $p$ is a ring epimorphism and since $p \circ i =\id_{R_{\Tilde{\Delta}}}$, the map $i$ is 
a bigraded algebra retract. We may apply a result from \cite{H} to the bigraded situation. 
It yields that  
$P^R_K= P^R_{R_{\Tilde{\Delta}}} P^{R_{\Tilde{\Delta}}}_K$ where  
$R_{\Tilde{\Delta}} = R/\nn_y$ is considered as a bigraded $R$-module. Since $R$ and 
$R_{\Tilde{\Delta}}$ are Koszul, the equality of bigraded Poincar\'e 
series implies that $\nn_y$ has a bigraded $1$-linear $R$-resolution. This concludes the
proof. 
\end{proof}

\begin{Proposition}\label{linear}
Let $c,d \geq 0$ be two integers. If $\nn_x$ and $\nn_y$ have bigraded linear resolutions, 
then
\begin{enumerate}
\item
the sidediagonal module $R^{(c,d)}_{\Delta}$ has a linear $R_{\Delta}$-resolution.
\item
the relative Veronese module $R^{(c,d)}_{\Tilde{\Delta}}$ has 
a bigraded linear $R_{\tilde{\Delta}}$-resolution.
\end{enumerate}
\end{Proposition}

For the proof of the proposition we need a fact which is stated in \cite{CHTV}. 
\begin{Lemma}\label{reg}
Let $A$ be a standard graded $K$-algebra, $M$ a finitely generated $A$-module and
$$\ldots \to N_r \to N_{r-1} \to \ldots \to N_1 \to N_0 \to M \to 0$$
be an exact complex of finitely generated graded $A$-modules. Then: 
\begin{enumerate}
\item
Let $h \in \NN$ and let $a \in \ZZ$ such that $t_s(N_r) \leq a+r+s$ 
for all $0 \leq r \leq h$ and $0 \leq s \leq h-r$. Then $t_h(M) \leq a+h$.

\item
$\reg_A M \leq \sup\{\reg_A N_i - i \colon i \in \NN\}$.
\end{enumerate}
\end{Lemma}

\begin{proof}[Proof of Proposition \rm \ref{linear}]
Since the proofs of (a) and (b) are similar, we only consider part (a). Moreover,
is enough to show that all modules $R^{(c,d)}_{\Delta}$ with $(c,d) \in \mathcal{I}(a,b)$
have linear resolutions. Let $G_\pnt$ be the bigraded minimal free $R$-resolution of 
$\nn_x$. Since $\nn_x$ has a bigraded $1$-linear resolution, every free module $G_r$ is of
the form  
$$G_r= \bigoplus_{p+q=r+1 \atop p \geq 1} R(-p,-q)^{\beta_{r,(p,q)}}$$
with non-negative integers $\beta_{r,(p,q)}$. Observe that for $c \geq 1$ 
and $(c,d) \in \mathcal{I}(a,b)$ we have $(\nn_x)^{(c,d)}=R^{(c,d)}_\Delta$. Applying the 
functor $-^{(c,d)}_\Delta$, we obtain an acyclic complex 
$(G_\pnt)^{(c,d)}_\Delta \to R^{(c,d)}_\Delta \to 0$ where
$$(G_r)^{(c,d)}_\Delta = \bigoplus_{p+q=r+1 \atop p \geq 1} (R^{(c_{p,q},d_{p,q})}_\Delta 
(-l_{p,q}))^{\beta_{r,(p,q)}}$$
By Remark \ref{shifts}(a), all occurring shifts $l_{p,q}$ are at most $r$. 
Similarly, let $H_\pnt$ be the minimal free resolution of $\nn_y$.
Then for $d \geq 1$ and $(c,d) \in \mathcal{I}(a,b)$ we observe that 
$(\nn_y)^{(c,d)}_\Delta= R^{(c,d)}_\Delta$ and the shifts in $(H_r)^{(c,d)}_\Delta$ are 
bounded by $r$.

To conclude the proof we show by induction that $t_h(R^{(c,d)}_\Delta) \leq h$
for all $h \in \NN$ and $(c,d) \in \mathcal{I}(a,b)$. First we use induction on $h$.
The modules $R^{(c,d)}_\Delta$ are generated in degree $0$, thus 
$t_0(R^{(c,d)}_\Delta) = 0$. Let now $h \geq 1$. We apply induction on $c+d$ where 
$(c,d) \in \mathcal{I}(a,b)$. For $c+d=0$ it follows that $c=0$ and $d=0$ and therefore
trivially $t_h( R_\Delta ) \leq h$. Let now $c+d > 0$. Then $c \geq 1$ or $d \geq 1$. 

We discuss the case $c \geq 1$ first. In order to apply Lemma \ref{reg}(a) to the exact 
complex $(G_\pnt)^{(c,d)}_\Delta \to R^{(c,d)}_\Delta \to 0$  we show that 
$t_s((G_r)^{(c,d)}_\Delta) \leq r+s$ for all $0 \leq r \leq h$ and $0 \leq s \leq h - r$.
Observe that $(G_0)^{(c,d)}_\Delta$ is a direct sum of $n$ copies of $(R^{(c-1,d)}_\Delta)$.    
Since $(c-1,d) \in \mathcal{I}(a,b)$, the induction hypothesis on $c+d$ implies that 
$t_s((G_0)^{(c,d)}_\Delta ) \leq s$ for all $0 \leq s \leq h$. 
For $1 \leq r \leq h$ and $0 \leq s \leq h - r$, we have 
$$t_s((G_r)^{(c,d)}_\Delta ) \leq t_s \Bigl(\bigoplus_{p+q=i+1 \atop p \geq 1} 
R^{(c_{p,q},d_{p,q})}_\Delta \Bigr) + r \leq s+r$$
where the first inequality holds because $l_{p,q} \leq r$ for all occurring $p,q$ and the 
second inequality holds by induction on $h$. Now Lemma \ref{reg} implies that 
$t_h(R^{(c,d)}_\Delta) \leq h$.

If $c=0$ and $d \geq 1$ the argument above similarly applies to the complex
$(H_\pnt)^{(c,d)}_\Delta \to R^{(c,d)}_\Delta \to 0$.
\end{proof}

As a direct consequence of Lemma \ref{Koszul_linear} and Lemma \ref{linear} we obtain

\begin{Corollary}\label{relative}
Let $c,d \geq 0$ be two integers. If $R$ is Koszul, then all sidediagonal modules 
$R^{(c,d)}_{\Delta}$ have linear $R_\Delta$-resolutions and all relative Veronese 
modules $R^{(c,d)}_{\Tilde{\Delta}}$ have bigraded linear $R_{\Tilde{\Delta}}$-resolutions.
\end{Corollary}

We use this corollary to get upper bounds for the regularity of sidediagonal
and relative Veronese modules. 

\begin{Theorem}\label{Koszul_reg}
Let $R$ be Koszul, $M$ be a finitely generated, bigraded $R$-module with $r=\reg_R M$
and $\indeg(M)=0$. 
\begin{enumerate}
\item
Let $(c,d) \in \mathcal{I}(a,b)$. Then
$$\reg_{R_\Delta} M_\Delta^{(c,d)} \leq \begin{cases}
				  				   		  \max\{ 0, \lceil \frac{r-c}{a} \rceil \},  
										  & \text{ if } b=0 \text{ and } a > 0, \cr
										  \max\{ 0, \lceil \frac{r-d}{b} \rceil \},  
										  & \text{ if } a=0 \text{ and } b > 0, \cr
										  \max\{ 0, \lceil \frac{r-c}{a} \rceil,
										  \lceil \frac{r-d}{b} \rceil \},
										  & \text{ if } a,b \geq 1. \cr 
				  				   		\end{cases} $$
\item
Let $(c,d) \in \Tilde{\mathcal{I}}(a,b)$. Then
$$\reg_{R_{\Tilde{\Delta}}} M_{\Tilde{\Delta}}^{(c,d)} \leq 
						   				 \begin{cases}
										    \max\{ 0, \lceil \frac{r}{a} \rceil \},  
										    & \text{ if } b=0 \text{ and } a > 0, \cr
										    \max\{ 0, \lceil \frac{r}{b} \rceil \},  
										    & \text{ if } a=0 \text{ and } b > 0, \cr

                                            \min\{r,\lceil \frac{r-c}{a} - \frac{d}{b} + 1 
											\rceil\},  & \text{ if } 1 \leq a \leq b. \cr 
				  				         \end{cases} $$
In particular, if $1 \leq a \leq b$, then 
$\reg_{R_{\Tilde{\Delta}}} M \leq \min\{r, \lceil \frac{r}{a}+1 \rceil \}$.
\end{enumerate}
\end{Theorem}
\begin{proof}
Let $F_\pnt$ be the minimal free $R$-resolution of $M$. 
Since $\reg_R M =r$ we have
$F_i =\bigoplus_{i \leq p+q \leq i+r}  R(-p,-q)^{\beta_{i,(p,q)}}$
for some non-negative integers $\beta_{i,(p,q)}$. 
For the proof of part (a) we restrict to the case $a, b \geq 1$. 
The other cases follow similarly. By Remark \ref{shifts}(a) and Corollary 
\ref{relative} we observe that

$\reg_{R_\Delta} (F_i)^{(c,d)}_\Delta \leq \max\{0,\lceil \frac{i+r-c}{a} \rceil,
   				   							      \lceil \frac{i+r-d}{b} \rceil\}
	\leq \max\{0,\lceil \frac{r-c}{a} \rceil, \lceil \frac{r-d}{b} \rceil\} + i$

\noindent
Now the claim follows from Lemma \ref{reg}(b). 
For part (b) we also restrict to the case $a \geq 1$ and $b \geq 1$.  
Use Remark \ref{shifts}(b) and Corollary \ref{relative} to observe that

$\reg_{R_{\Tilde{\Delta}}} (F_i)^{(c,d)}_{\Tilde{\Delta}} \leq 
                   \max \{ \max\{0,\lceil \frac{p-c}{a} \rceil \}
   				  + \max\{0,\lceil \frac{q-d}{b} \rceil\} \colon i \leq p+q \leq i+r \}$
									  
\noindent
The claim follows from a case by case computation using $1 \leq a \leq b$ and 
Lemma \ref{reg}(b). Then the upper bound for $\reg_{R_{\Tilde{\Delta}}} M$ follows from 
the fact that $M$ decomposes into the finite sum
$M =\bigoplus_{(c,d) \in \Tilde{I}(a,b)} M^{(c,d)}_{\Tilde{\Delta}}.$
\end{proof}

As a direct consequence of Theorem \ref{Koszul_reg} the modules $M_\Delta$ and 
$M_{\Tilde{\Delta}}$ have small regularities for $a,b \gg 0$. More concrete, we have

\begin{Corollary}\label{large_reg}
Let $M$ be a finitely generated, bigraded $R$-module.
\begin{enumerate}
\item
If $\max\{a,b\} \geq \reg_R M$, then $\reg_{R_\Delta} M_\Delta \leq \min\{1, \reg_R M \}$. 
\item
Let $a,b \geq 1$. If $\min\{a,b\} \geq \reg_R M$, then 
$\reg_{R_{\Tilde{\Delta}}} M_{\Tilde{\Delta}} \leq 
\min\{2, \reg_R M\}$ and
$\reg_{R_{\Tilde{\Delta}}} M \leq \min\{2, \reg_R M\}$. 
\end{enumerate}
\end{Corollary}  

The main result follows immediately from the results above.

\begin{proof}[Proof of Theorem \rm \ref{main}]
Note that a graded $K$-algebra $A$ is Koszul if and only if $\reg_A K = 0$. 
Since $K_\Delta = K$ and $K_{\Tilde{\Delta}}=K$, the claim follows from Theorem
\ref{Koszul_reg}.
\end{proof}

Note that the converse of Theorem \ref{main} is false. Take, for example, 
the algebra $R= K[x_1,y_1]/(x_1 y_1^2)$. Since the defining ideal of $R$ is generated in
degree $3$, $R$ is not Koszul. But every diagonal $R_\Delta$ is Koszul because 
$R_\Delta$ is either isomorphic to a polynomial ring 
$K[t]$, to $K$ or to the Koszul algebra $K[t]/(t^2)$.

\section{Applications}
In this section we present some applications which arise naturally in the study of 
symmetric algebras and Rees algebras.  
In the following $A$ will always denote a positively graded algebra, i.e. 
$A = K[x_1,\ldots,x_n]/Q$ where $\deg(x_i)=1$ and $Q \subset K[x_1,\ldots,x_n]$ is 
a homogeneous ideal. Let $\mm \subset A$ be the graded maximal ideal.   

We first consider symmetric algebras. Let $M$ be a graded $A$-module with homogeneous 
generators $f_1,\ldots,f_m$ and let $(a_{ij})_{i=1,\ldots,t \atop j=1,\ldots,m}$ be the 
corresponding relation matrix. The symmetric algebra $S(M) = \bigoplus_{j \geq 0} S^j(M)$ of 
$M$ has a presentation of the form $S(M)=A[y_1,\ldots,y_m]/J$ where 
$J=(g_1, \ldots,g_t)$ and $g_i = \sum_{j=1}^{m} a_{ij} y_j$ for $i=1,\ldots,t$. If  
$f_1,\ldots,f_m$  have the same degree, then $S(M)$ is standard bigraded by assigning the
degree $(1,0)$ to the residue class of $x_i$ for $i=1,\ldots,n$ and setting 
$\deg(y_i)=(0,1)$. Note that $S^j(M)$ is a graded $A$-module. As an application of the 
main result we obtain

\begin{Corollary}\label{Koszul_sym}
If $S(M)$ is Koszul, then $A$ is Koszul and the module $S^j(M)$ has a linear 
resolution for all $j \geq 0$.
\end{Corollary}

\begin{proof}
Let $\Delta$ be the $(1,0)$-diagonal. Then $S(M)_\Delta = A$, and  
$S^j(M)=S(M)^{(0,j)}_\Delta$. Thus $A$ is Koszul and, by Corollary \ref{relative}, 
the module $S^j(M)$ has a linear $A$-resolution.
\end{proof}

As one might expect it seems to be a strong condition that the symmetric algebra $S(M)$ is 
Koszul. In a more specific case however, when $M=\mm$ is the graded maximal ideal of 
a Koszul algebra, we have a sufficient condition.     

\begin{Theorem}\label{golod}
Let $A = K[x_1,\ldots,x_n] / Q$ and $K$ be an infinite field with $\ch(K) \neq 2$. 
If $Q$ has a $2$-linear resolution over $K[x_1,\ldots,x_n]$, then the defining ideal of 
$S(\mm)$ has a quadratic Gr\"obner basis with respect to a reverse lexicographic
term order. In particular, $S(\mm)$ is Koszul.
\end{Theorem}

It is well-known that the existence of a quadratic Gr\"obner basis for the defining
ideal of an algebra implies the Koszul property. For details on Gr\"obner bases and 
generic initial ideals refer to \cite{E}.

\begin{Lemma}\label{gin}
Let $Q \subset K[x_1,\ldots,x_n]$ be an homogeneous ideal, $K$ infinite and $\ch(K) \neq 2$.
Moreover, let $\Gin(Q)$ denote the generic initial ideal with respect to the reverse 
lexicographic term order induced by $x_1 > \ldots > x_n$. If $Q$ has a $2$-linear resolution, 
then $\Gin(Q)$ is quadratic and satisfies the following condition:

\noindent
$(*) \quad$ If $x_i x_j \in \Gin(Q)$ and $i \leq j$, then $x_i x_k \in \Gin(Q)$ for all 
$k < j$.  
\end{Lemma}

\begin{proof}
We use some results about $\Gin(Q)$. It is known that $\Gin(Q)$ is a Borel-fixed ideal and
that $\reg\Gin(Q)=\reg(Q)=2$ (see \cite[$20.21$]{E}). Thus $\Gin(Q)$ is quadratic. Since 
$\ch(K) \neq 2$ we obtain that ${\Gin(Q)}_2$ is stable (see \cite[$15.23$b]{E}) and 
therefore satisfies condition $(*)$.
\end{proof} 

We need some notation taken from \cite{HPT}.
Let $S=K[x_1,\ldots,x_n,y_1,\ldots,y_n]$ be the standard bigraded polynomial ring and
$f = \sum_{1 \leq i_1 \leq i_2 \leq \ldots \leq i_d \leq n} a_{i_1 i_2 \cdots i_d} 
x_{i_1} x_{i_2} \cdots x_{i_d}$ a form of degree $(d,0)$. We set
$$f^{(k)}=  \sum_{1 \leq i_1 \leq i_2 \leq \ldots \leq i_d \leq n} a_{i_1 i_2 \cdots i_d} 
x_{i_1} x_{i_2} \cdots x_{i_{d-k}} y_{i_{d-k+1}}\cdots y_{i_d}$$
for $k=0,\ldots,d$. Note that $f^{(k)}$ is bihomogeneous of degree $(d-k,k)$.
Moreover, let $\delta_{ij} = x_i y_j - x_j y_i$ for $i \neq j$ and
$L=\{\delta_{ij} \colon i \neq j\}$. We need the following lemma.

\begin{Lemma}\label{ini_lemma}
Let $>$ denote the reverse lexicographic term order on $S$ induced by 
$x_1 > x_2 > \ldots > x_n > y_1 > \ldots > y_n$ and $\phi : S \to S$ be the 
homomorphism with $\phi(x_i) = x_i$ and $\phi(y_i) =x_i$ for $i=1,\ldots,n$.
Assume that $f \in S$ is a bihomogeneous polynomial of degree $(s,t)$ such that  
$\ini(f) = x_{i_1} x_{i_2} \cdots x_{i_s} y_{j_1} \cdots y_{j_t}$ satisfies
$i_1 \leq i_2 \leq \ldots \leq i_s \leq j_1 \leq \ldots \leq j_t$. Then 
$\ini(\phi(f)) = \phi(\ini(f))$ and $\ini(f) = \phi (\ini(f))^{(t)}$.
\end{Lemma}

\begin{proof}
With the condition $i_1 \leq i_2 \leq \ldots \leq i_s \leq j_1 \leq \ldots \leq j_t$
it is easy to see that $\phi ( \ini(f)) > \phi (v)$ for all monomials $v$ of $f$ with 
$v < \ini(f)$.
\end{proof}

\begin{proof}[Proof of Theorem \rm\ref{golod}]
By Lemma \ref{gin}, we may assume that the defining ideal $Q$ of $A$ has a quadratic 
Gr\"obner basis $g_1,\ldots,g_t$ with respect to the reverse lexicographic term order
induced by $x_1 > x_2 > \ldots > x_n$ and that $\ini(Q)$ satisfies condition \ref{gin} $(*)$.
It is easy to see that $S(\mm)$ has a presentation $S(\mm) = S/J$ where 
$J=(g_1,\ldots,g_t, g_1^{(1)},\ldots,g_t^{(1)},L)$. 

Let $>$ denote the reverse lexicographic term order on $S$ induced by 
$x_1 > x_2 > \ldots > x_n > y_1 > \ldots > y_n$. We will show that the set
$G =\{ g_1,\ldots,g_t, g_1^{(1)},\ldots,g_t^{(1)}\} \union L$
is a Gr\"obner basis for $J$ with respect to $>$ which concludes the proof of the 
theorem.

Let $f \in J$ be a bihomogeneous polynomial of degree $(s,t)$. Then $s \geq 1$. 
We show that $\ini(f)$ is divided by some $\ini(g)$ with $g \in G$.
Let $\ini(f) = x_{i_1} x_{i_2} \cdots x_{i_{s}} y_{j_1} \cdots y_{j_t}$
where $i_1 \leq i_2 \leq \ldots \leq i_s$ and $j_1 \leq j_2 \leq \ldots \leq j_t$
If there exist indices $p,q$ such that $i_p > j_q$ then $\ini (\delta_{i_p j_q})$
divides $\ini(f)$ which is the claim.

Otherwise we have 
$i_1 \leq i_2 \leq \ldots \leq i_s \leq j_1 \leq j_2 \leq \ldots \leq j_t$.
Let $\phi$ denote the homomorphism from Lemma \ref{ini_lemma}.
Since $f \in J$, it follows that $\phi(f) \in Q$. By Lemma \ref{ini_lemma}, we have 
$\ini(f) = \ini(\phi(f))^{(t)}$ where 
$\ini(\phi(f)) = x_{i_1} x_{i_2} \cdots x_{i_{s}} x_{j_1} \cdots x_{j_t}$. Since 
$g_1,\ldots,g_t$ is a Gr\"obner basis for $Q$, there exists a 
$g \in \{g_1,\ldots,g_t\}$ such that $\ini(g)$ divides $\ini(\phi(f))$. 
By condition \ref{gin} $(*)$, we may assume that 
$\ini(g) = x_{i_1} x_{i_2}$, if $s>1$, or
$\ini(g) = x_{i_1} x_{j_1}$, if $s=1$. Now $\ini(g)$ or $\ini(g^{(1)})$ divides 
$\ini(f)$.
\end{proof} 

Under the strong assumption of Theorem \ref{golod} it follows from Corollary 
\ref{Koszul_sym} that $S^j(\mm)$ has linear resolution for all $j \geq 1$. Actually we 
have

\begin{Proposition}\label{m_symmetric}
\rm Let $j \geq 1$. If $A$ is Koszul, then $S^{j}(\mm)$ has a linear $A$-resolution.
\end{Proposition}

In the proof we use results from \cite{HSV} and
some basic facts about the Koszul complex (see \cite[Section I.1.6]{BH} for details).
   
\begin{proof}
Let $S_x=K[x_1,\ldots,x_n]$ the standard graded polynomial ring and $A=S_x/Q$.
We may assume that the defining ideal $Q$ of $A$ does not contain linear forms.
Then $Q$ is generated in degree $2$. Let $\mm=(x_1,\ldots,x_n) \subset A$ be the graded 
maximal ideal of $A$. We denote by $K_\pnt$ the Koszul complex of the sequence 
$x_1,\ldots,x_n \in A$. Let $H_1(K_\pnt)$ be the first
homology group of this complex.
Recall that $S(\mm)=A[y_1,\ldots,y_n]/J$ for some bihomogeneous ideal $J$
and that $S^j(\mm)$ is generated by the residue classes of all monomials in degree $(0,j)$. 
We consider $S^j(\mm)$ as an $A$-module generated in degree $j$. For $j \geq 1$, there 
exists the downgrading homomorphism $\alpha_j : S^j(\mm) \to \mm S^{j-1}(\mm)$ mapping
a residue class of $y_{i_1} y_{i_2} \ldots y_{i_j}$ to the residue class of 
$x_{i_1} y_{i_2} \ldots y_{i_j}$ (see \cite[Section 2]{HSV}). 
Note that it does not matter 
which of the factors $y_{i_l}$ is replaced by $x_{i_l}$.  
 
To show that $S^j(\mm)$ for all $j \geq 1$ we use induction on $j$. For $j=1$, we have 
$S^1(\mm)= \mm$ which has a linear resolution because $A$ is Koszul. Let now $j>1$. 
We have the short exact sequence
\begin{equation}\label{downgrading}
0 \to U \to S^{j}(\mm) \overset{\alpha_j}{\to} \mm S^{j-1}(\mm) \to 0
\end{equation} 
where $U=\ker \alpha_j$. By \cite[Lemma 2.2]{HSV}, $U$ is a subquotient of the module 
$N = H_1(K_\pnt) \otimes_{A/\mm} (A/\mm)(-j+2)^s$ for some integer $s \geq 1$. The module $N$ is 
annihilated by $\mm$. Since $Q$ is generated in degree $2$,
it follows that $H_1(K_\pnt) \cong \Tor^{S_x}_1(A,K)$ is generated in degree $2$.
Therefore, $U \cong K(-j)^{t}$ for some integer $t \geq 0$ and $U$ has a $j$-linear
$A$-resolution because $A$ is Koszul. By induction hypothesis, $S^{j-1}(\mm)$ 
has $(j-1)$-linear $A$-resolution. Thus, by \cite[Lemma 6.4]{CHTV}, the module 
$\mm S^{j-1}(\mm)$ has a $j$-linear $A$-resolution. The assertion follows when we apply the
long exact sequence of the functor $\Tor^A(\cdot,K)$ to the sequence $(\ref{downgrading})$. 
\end{proof}

The hypothesis of Theorem \ref{golod} cannot be weakened to the assumption that
$A$ is only Koszul. A counter example is the algebra $A=K[x_1,x_2]/(x_1^2,x_2^2)$.
As a complete intersection, $A$ is Koszul, but with the help of the program Macaulay2 we
find that $S(\mm)$ is not Koszul. This example shows also that
the converse of Corollary \ref{Koszul_sym} is false because, by Proposition 
\ref{m_symmetric}, all symmetric powers $S^j(\mm)$ have a linear resolution.

A further intensively studied class of bigraded algebras are Rees algebras.
Let $I \subset A$ be a homogeneous ideal which is minimally generated by homogeneous 
elements $f_1, f_2 \ldots, f_m$ of the same degree $d$. Recall that the Rees ring 
$R(I) = A[It]$ of $I$ admits a standard bigrading assigning 
the degree $(1,0)$ to the generators of $\mm \subset A$  and setting 
$\deg(f_i t)=(0,1)$ for $i =1,\ldots,m$. As a consequence of Theorem \ref{main} we
observe  

\begin{Corollary}
If $R(I)$ is Koszul, then $A$ is Koszul and the ideal $I^j$ has a linear 
$A$-resolution for all $j \geq 0$. 
\end{Corollary}

\begin{proof}
Let $\Delta$ be the $(1,0)$-diagonal. Then $R(I)_\Delta = A$ and 
$I^j = R(I)^{(0,j)}_\Delta(-dj).$ Thus, by Corollary \ref{Koszul_linear}, the ideal 
$I^j$ has a linear $A$ resolution.
\end{proof}

Note that the Rees algebra $R(\mm)$ is always Koszul because it is a Segre product of
two Koszul algebras. In \cite{HPT} Herzog, Popescu and Trung prove that the 
defining ideal of $R(\mm)$ has a quadratic Gr\"obner basis provided $Q$ has a quadratic
Gr\"obner basis.
 
As direct consequences of Theorem \ref{main} we get some well-known facts about positively 
graded Koszul algebras. Let $d \geq 1$ be an integer. Recall that the $d^{th}$ Veronese 
subring $A^{(d)}$ of $A$ is the positively graded algebra 
$A^{(d)}= \bigoplus_{i \geq 0} A_{id}.$
For two standard graded $K$-algebras $A$ and $B$ the tensor product 
$A \otimes_K B = \bigoplus_{i,j \geq 0}  A_i \otimes_K B_j$
is a standard bigraded algebra. The Segre product of $A$ and $B$, denoted with $A \star B$, 
is the $(1,1)$-diagonal of $A \otimes_K B$. We recover some well-known results 
(see \cite{F}).

\begin{Corollary}
Tensor products, Segre products and Veronese subrings of Koszul algebras are Koszul.
\end{Corollary}

\begin{proof}
Let $F_\pnt$ (resp.\ $G_\pnt$) be the minimal free resolution of $K$ over $A$ (resp.\ $B$).
Then the tensor product $G_\pnt \otimes_K F_\pnt$ gives a minimal free resolution of $K$
over $A \otimes_K B$. Thus if $A$ and $B$ are Koszul, $A \otimes_K B$ is Koszul.
Now the Segre product $A \star B$ is the $(1,1)$-diagonal $A \otimes_K B$ which is Koszul
by Theorem \ref{main}.

Let $A$ be a positively graded Koszul algebra. Consider $A$ as a standard bigraded algebra 
where all generators have degree $(1,0)$. Then $A^{(d)}$ is a diagonal of $A$ and, by 
Theorem \ref{main}, $A^{(d)}$ is Koszul. 
\end{proof}

Let $M$ be a graded $A$-module. Recall from \cite{ABH} that the rate of $M$ is given by
$ \rate_A M = \sup\{ \frac{t_i(M)}{i} \colon i \geq 0 \}$.
A similar definition is given in \cite{B} where Backelin  proves that 
$A^{(d)}$ is Koszul  for $d \gg 0$.  Note an $A$-module $M$ is naturally
an $A^{(d)}$-module. Aramova, Barcanescu and Herzog have proved in \cite{ABH}  that 
\begin{equation}\label{rate}
\rate_{A^{(d)}} M \leq \lceil \rate_A M / d \rceil
\end{equation}
for an arbitrary $K$-algebra $A$ and all $d \geq c$ where $c$ is a constant depending
on $A$. Moreover, they showed that $c=1$, if $A$ is a polynomial ring. For this, they
used that the relative Veronese modules $A^{(d,j)} = \bigoplus_{i \geq 0} A_{id+j}$
for $j=0,\ldots,d-1$ have linear $A$-resolutions.
Since the relative Veronese modules coincide with sidediagonal modules, 
it follows from Corollary \ref{linear} that (\ref{rate}) is valid for $c \geq 1$ 
provided $A$ is Koszul. We get similar upper bounds for the regularity over Koszul 
algebras.
 
\begin{Corollary}
Let $A$ be Koszul and $M$ be a finitely generated graded $A$-module. Then
$\reg_{A^{(d)}} M \leq \lceil \reg_A M / d \rceil$ for all $d \geq 1$. In particular, 
$\reg_{A^{(d)}} M \leq 1$ if $d \geq \reg_A M$. 
\end{Corollary}

\begin{proof}
Consider $A$ as a bigraded algebra generated in degree $(1,0)$. Let $\Delta$ be the
$(d,0)$-diagonal of $A$. Then $A^{(d)}=A_\Delta$ and, as an $A^{(d)}$-module, we have
$M=\bigoplus_{c=0}^{d-1} M^{(c,0)}_\Delta$. By Theorem \ref{Koszul_reg}, the claim 
follows. 
\end{proof}

\section{Semigroup rings}
Finally, we study the consequences of the main result for bihomogeneous semigroup rings.
Let $\Lambda \subset \NN^d$ be a finitely generated semigroup. We call $\Lambda$
standard bigraded if 
\begin{enumerate}
\item
$\Lambda$ is the disjoint union $\bigcup_{i,j \geq 0} \Lambda_{(i,j)}$ and
\item 
$\Lambda_{(0,0)} =0$, $\Lambda_{(i,j)} + \Lambda_{(k,l)} \subset \Lambda_{(i+k,j+l)}$ for
all integers $i,j,k,l \geq 0$ and 
\item
$\Lambda$ is generated by elements of $\Lambda_{(1,0)}$ and $\Lambda_{(0,1)}$. 
\end{enumerate}  
We call the elements of $\Lambda_{(i,j)}$ bihomogeneous of degree $(i,j)$.
Similarly, one defines a graded semigroup.
Let $\Lambda$ be a standard bigraded semigroup which is minimally generated by 
$\alpha_1,\ldots,\alpha_n \in \Lambda_{(1,0)}$ and 
$\beta_1,\ldots,\beta_m \in \Lambda_{(0,1)}$ and let $K[t_1,\ldots,t_d]$ denote 
the polynomial ring. For every semigroup element 
$\lambda = (a_1,\ldots,a_d) \Lambda$ we associate the monomial 
$t^\lambda = t_1^{a_1} t_2^{a_2} \cdots t_d^{a_d}$. 
Recall that the semigroup ring $K[\Lambda]$ is the $K$-algebra 
generated by the monomials $t^{\alpha_i}, t^{\beta_j}$ where $i=1,\ldots,n$
and $j=1,\ldots,m$. Let $\phi : S \to K[\Lambda]$
be the epimorphism with $\phi(x_i) = t^{\alpha_i}$ and $\phi(y_j) = t^{\beta_j}$.
Then $J = \ker(\phi)$ is called the toric ideal of the semigroup ring $K[\Lambda]$. If 
$\Lambda$ is bigraded then $K[\Lambda] = S/J$ is a standard bigraded algebra.  

The divisibility relation of the monomials in $K[\Lambda]$ defines a partial
order $\preceq$ on $\Lambda$. For $\mu, \lambda \in \Lambda$ we set $\mu \preceq \lambda$
if $\lambda = \sigma + \mu$ for some $\sigma \in \Lambda$.
Then the open intervals
$(\mu ,\lambda) = \{ \sigma \in \Lambda \colon \mu \prec \sigma \prec \lambda \}$
are partially ordered with the induced ordering.

Let $(P, \preceq)$ be a finite poset. Recall that the boundary complex $\Gamma(P)$ is
the simplicial complex whose faces are the totally ordered subsets of $P$.
For $\lambda \in \Lambda$ we denote the boundary complex of the interval 
$(0,\lambda)$ by  $\Gamma_\lambda$. The following is stated in \cite[Corollary 2.2]{PRS} 
and \cite{HRW}.

\begin{Proposition}
$K[\Lambda]$ is Koszul if and only if $\Gamma_\lambda$ is Cohen-Macaulay for all
$\lambda \in \Lambda$. 
\end{Proposition}  

Let $\Lambda$ be a bigraded semigroup. In analogy to the definition for $K$-algebras we set 
$$ \Lambda_\Delta = \bigcup_{i \geq 0} \Lambda_{(ia,ib)} \quad \text{and} \quad
\Lambda_{\Tilde{\Delta}} = \bigcup_{i,j \geq 0} \Lambda_{(ia,jb)}$$
for the $(a,b)$-diagonal $\Delta$. Note that $\Lambda_\Delta$ is graded and partially
ordered by the induced ordering. If $\lambda \in \Lambda_{(ia,ib)}$, then  we use 
$(\Gamma_\lambda)_\Delta$ for the boundary complex of the induced open interval 
$(0,\lambda) \subset \Lambda_\Delta$. Similarly, we define 
$(\Gamma_\lambda)_{\Tilde{\Delta}}$ for $\lambda \in \Lambda_{(ia,jb)}$.
Finally, we reformulate our main result for semigroup rings.

\begin{Corollary}
Let $\Lambda \subset \NN^d$ a bigraded semigroup and $\Delta$ a diagonal. If 
$\Gamma_\lambda$ is Cohen-Macaulay for all $\lambda \in \Lambda$, then
\begin{enumerate}

\item
$(\Gamma_\lambda)_\Delta$ is Cohen-Macaulay for all $\lambda \in \Lambda_\Delta$.

\item
$(\Gamma_\lambda)_{\Tilde{\Delta}}$ is Cohen-Macaulay for all 
$\lambda \in \Lambda_{\Tilde{\Delta}}$.

\end{enumerate}
\end{Corollary}

\bigskip

\noindent
Stefan Blum\\
FB6 Mathematik und Informatik\\
Universit\"at-GH Essen\\
Postfach 103764\\
45117 Essen\\
Germany\\
stefan.blum@uni-essen.de
\end{document}